\documentclass[12pt]{article}
\usepackage{amsmath,amsfonts,amssymb,amsthm}
\input amssym.def
\begin{document}
\def \Z{\Bbb Z}
\def \C{\Bbb C}
\def \R{\Bbb R}
\def \Q{\Bbb Q}
\def \N{\Bbb N}

\def \A{{\mathcal{A}}}
\def \D{{\mathcal{D}}}
\def \E{{\mathcal{E}}}
\def \H{\mathcal{H}}
\def \S{{\mathcal{S}}}
\def \V{{\mathcal{V}}ir}
\def \wt{{\rm wt}}
\def \tr{{\rm tr}}
\def \span{{\rm span}}
\def \Res{{\rm Res}}
\def \Der{{\rm Der}}
\def \End{{\rm End}}
\def \Ind {{\rm Ind}}
\def \Irr {{\rm Irr}}
\def \Aut{{\rm Aut}}
\def \GL{{\rm GL}}
\def \Hom{{\rm Hom}}
\def \mod{{\rm mod}}
\def \ann{{\rm Ann}}
\def \ad{{\rm ad}}
\def \rank{{\rm rank}\;}
\def \<{\langle}
\def \>{\rangle}

\def \g{{\frak{g}}}
\def \h{{\hbar}}
\def \k{{\frak{k}}}
\def \sl{{\frak{sl}}}
\def \gl{{\frak{gl}}}

\def \be{\begin{equation}\label}
\def \ee{\end{equation}}
\def \bex{\begin{example}\label}
\def \eex{\end{example}}
\def \bl{\begin{lem}\label}
\def \el{\end{lem}}
\def \bt{\begin{thm}\label}
\def \et{\end{thm}}
\def \bp{\begin{prop}\label}
\def \ep{\end{prop}}
\def \br{\begin{rem}\label}
\def \er{\end{rem}}
\def \bc{\begin{coro}\label}
\def \ec{\end{coro}}
\def \bd{\begin{de}\label}
\def \ed{\end{de}}

\newcommand{\m}{\bf m}
\newcommand{\n}{\bf n}
\newcommand{\nno}{\nonumber}
\newcommand{\nord}{\mbox{\scriptsize ${\circ\atop\circ}$}}
\newtheorem{thm}{Theorem}[section]
\newtheorem{prop}[thm]{Proposition}
\newtheorem{coro}[thm]{Corollary}
\newtheorem{conj}[thm]{Conjecture}
\newtheorem{example}[thm]{Example}
\newtheorem{lem}[thm]{Lemma}
\newtheorem{rem}[thm]{Remark}
\newtheorem{de}[thm]{Definition}
\newtheorem{hy}[thm]{Hypothesis}
\makeatletter \@addtoreset{equation}{section}
\def\theequation{\thesection.\arabic{equation}}
\makeatother \makeatletter

\begin{center}
{\Large \bf Associating quantum vertex algebras to certain deformed
Heisenberg Lie algebras}
\end{center}

\begin{center}
{Haisheng Li\footnote{Partially supported by NSA grant
H98230-11-1-0161}\\
Department of Mathematical Sciences\\ Rutgers University, Camden, NJ
08102}
\end{center}

\begin{abstract}
We associate quantum vertex algebras and their $\phi$-coordinated
quasi modules to certain deformed Heisenberg algebras.
\end{abstract}

\section{Introduction}
It is now well known (see \cite{fz}; cf. \cite{ll}) that certain
infinite-dimensional Lie algebras such as affine (Kac-Moody) Lie
algebras, the Virasoro algebra, and Heisenberg algebras of a certain
type can be canonically associated with vertex algebras and their
modules. This association can be briefly outlined as follows: First,
the canonical generating functions of their generators are {\em
mutually local} in the sense that for any two generating functions
$a(x)$ and $b(x)$, there exists a nonnegative integer $k$ such that
\begin{eqnarray}\label{locality-intro}
(x_{1}-x_{2})^{k}a(x_{1})b(x_{2})=(x_{1}-x_{2})^{k}b(x_{2})a(x_{1}).
\end{eqnarray}
Second, a conceptual result of \cite{li-local} states that for any
vector space $W$, any local subset of $\Hom (W,W((x)))$,
alternatively denoted by $\E(W)$, generates a vertex algebra
canonically with $W$ as a faithful module. Moreover, the generated
vertex algebras have a canonical module structure identified as a
so-called vacuum module for the corresponding Lie algebra.

In this paper, we study two deformed Heisenberg Lie algebras, which
have appeared in the study of quantum algebras, and our goal is to
associate vertex algebras or their likes to these Lie algebras.
 The first we are
concerned about is the Heisenberg algebra (see \cite{efr}) with
generators $a_{n}$ for $n\in \Z$, subject to relations
\begin{eqnarray}
[a_{m},a_{n}]=[m]_{q}\delta_{m+n,0}
\end{eqnarray}
for $m,n\in \Z$, where the complex parameter $q$ is neither zero nor
a root of unity and where $[m]_{q}=\frac{q^{m}-q^{-m}}{q-q^{-1}}$,
while the second is the Heisenberg algebra (cf. \cite{bp}) with
generators $b_{n}$ for $n\in \Z$, subject to relations
\begin{eqnarray}
[b_{m},b_{n}]=m(1-q^{|m|})\delta_{m+n,0}.
\end{eqnarray}
Heisenberg Lie algebras are the simplest nonabelian Lie algebras,
however, as we shall see, associating vertex algebras or their likes
to them is by no means straightforward. The main purpose of this
paper is to provide a simple but illuminating example.

For the first one, forming a generating function $a(x)=\sum_{n\in
\Z}a_{n}x^{-n}$, we have
\begin{eqnarray}
[a(x_{1}),a(x_{2})]=\frac{1}{q-q^{-1}}\left(\delta\left(\frac{qx_{2}}{x_{1}}\right)
-\delta\left(\frac{x_{2}}{qx_{1}}\right)\right),
\end{eqnarray}
while for the second one, using generating function $b(x)=\sum_{n\in
\Z}b_{n}x^{-n}$ we have
\begin{eqnarray}
[b(x_{1}),b(x_{2})]=\left(x_{2}\frac{\partial}{\partial
x_{2}}\right)\delta\left(\frac{x_{2}}{x_{1}}\right)
+\frac{qx_{1}/x_{2}}{(1-qx_{1}/x_{2})^{2}}
-\frac{qx_{2}/x_{1}}{(1-qx_{2}/x_{1})^{2}}.
\end{eqnarray}
Notice that neither generating function $a(x)$ nor $b(x)$ forms a
local set. Thus the results of \cite{li-local} are {\em not}
applicable to these two Lie algebras.

In the past, for various purposes we have significantly generalized
the results of \cite{li-local} in several directions. One of the
generalizations was obtained in \cite{li-gamma} with a purpose to
associate vertex algebras to a family of infinite-dimensional Lie
algebras called quantum tori Lie algebras (see \cite{gkl},
\cite{gkk}). For these Lie algebras, their generating functions
satisfy a generalized locality relation
\begin{eqnarray}\label{quasi-local-intro}
p(x_{1},x_{2})a(x_{1})b(x_{2})=p(x_{1},x_{2})b(x_{2})a(x_{1})
\end{eqnarray}
for some nonzero homogeneous polynomial $p(x_{1},x_{2})$. Motivated
by this, we then introduced a notion of quasi locality, using
commutativity relation (\ref{quasi-local-intro}).  It was proved
that for any vector space $W$, every quasi local subset of $\E(W)$
generates a vertex algebra in a certain {\em new} way, with $W$ as
what we called a quasi module. In this way, we obtained a new
construction of vertex algebras together with a theory of quasi
modules for vertex algebras, not to mention that vertex algebras
were associated to quantum tori Lie algebras as an application.

As with the two aforementioned Heisenberg algebras, we see that the
generating function $a(x)$ of the first one forms a quasi local set,
so that one can use the results of \cite{li-gamma} to associate a
vertex algebra to this Lie algebra. However, the structure of the
associated vertex algebra is not as neat as expected. As for the
second one, since the generating function $b(x)$ does not give rise
to a quasi local set, the construction of \cite{li-gamma} is not
applicable. In fact, the notion of vertex algebra is not general
enough; we need a generalization of the notion of vertex algebra.

Inspired by Etingof-Kazhdan's notion of quantum vertex operator
algebra (see \cite{ek}), in \cite{li-qva1} we introduced a notion of
weak quantum vertex algebra and a notion of quantum vertex algebra,
greatly generalizing the notion of vertex superalgebra. Furthermore,
we established a conceptual construction. For a vector space $W$, a
subset $U$ of $\E(W)$ is said to be {\em $\S$-local} if for any
$u(x),v(x)\in U$, there exist
$$u^{(i)}(x),v^{(i)}(x)\in U,\ f_{i}(x)\in \C((x))\ \ (i=1,\dots,r)$$ such that
$$(x_{1}-x_{2})^{k}u(x_{1})v(x_{2})
=(x_{1}-x_{2})^{k}\sum_{i=1}^{r}f_{i}(x_{2}-x_{1})v^{(i)}(x_{2})u^{(i)}(x_{1})$$
for some nonnegative integer $k$. It was proved that any $\S$-local
subset of $\E(W)$ generates a weak quantum vertex algebra with $W$
as a faithful module.

Note that the generating functions of the two aforementioned
Heisenberg algebras do {\em not} give rise to $\S$-local sets. Thus
the results of \cite{li-qva1} are not applicable {\em either.} To
deal with the two Heisenberg algebras, especially the second one, we
shall need a different idea.

In the general field of vertex algebras, a fundamental problem (see
\cite{fj}, \cite{efk}) has been to develop a theory of quantum
vertex algebras so that quantum vertex algebras can be canonically
associated to quantum affine algebras. For many years, solving this
very problem has been the main driving force for us to study quantum
vertex algebras. Finally, we made a significant progress in
\cite{li-phi-module}. As the crucial step, we developed a theory of
what we called $\phi$-coordinated quasi modules for weak quantum
vertex algebras (actually for more general nonlocal vertex
algebras), where $\phi=\phi(x,z)\in \C((x))[[z]]$ is a formal series
satisfying
$$\phi(x,0)=x,\ \ \ \ \phi(\phi(x,z_{1}),z_{2})=\phi(x,z_{1}+z_{2}).$$
All such $\phi$ were completely determined therein, and two
particular examples are $\phi(x,z)=x+z$ and $\phi(x,z)=xe^{z}$. For
$\phi(x,z)=x+z$, the notion of $\phi$-coordinated quasi module
reduces to the notion of quasi module. What is important for us to
deal with quantum affine algebras is the case with $\phi=xe^{z}$. In
this case, we obtained a general construction. Let $W$ be a general
vector space. A subset $U$ of $\E(W)$ is said to be {\em quasi
$\S_{trig}$-local} if for any $a(x),b(x)\in U$, there exist
$$u^{(i)}(x),v^{(i)}(x)\in U,\ f_{i}(x)\in \C(x)\ \ (i=1,\dots,r),$$
where $\C(x)$ denotes the field of rational functions, such that
$$p(x_{1},x_{2})u(x_{1})v(x_{2})
=p(x_{1},x_{2})\sum_{i=1}^{r}f_{i}(x_{1}/x_{2})v^{(i)}(x_{2})u^{(i)}(x_{1})$$
for some nonzero polynomial $p(x_{1},x_{2})$. It was proved that any
quasi $\S_{trig}$-local subset of $\E(W)$ generates a weak quantum
vertex algebra in a certain way with $W$ as a $\phi$-coordinated
quasi module where $\phi(x,z)=xe^{z}$. As an application, we have
successfully associated weak quantum vertex algebras to quantum
affine algebras.

In this current paper, by applying the results of
\cite{li-phi-module} we associate vertex algebras and
$\phi$-coordinated quasi modules to the first Lie algebra, whereas
we associate quantum vertex algebras to the second. Note that by
\cite{li-phi-module} {\em weak} quantum vertex algebras can be
associated to them {\em conceptually}. The main task here is to
construct the desired weak quantum vertex algebras concretely and
show that they are indeed quantum vertex algebras. To explicitly
determine the associated vertex algebras and quantum vertex
algebras, we introduce and employ certain Heisenberg algebras. We
show that the associated quantum vertex algebras are simple. This as
an illustrating example shows how one can associate quantum vertex
algebras to more general infinite-dimensional Lie algebras.

This paper is organized as follows: In Section 2, we recall some
necessary notions and results. In Section 3, we study the first
deformed Heisenberg algebra in the context of vertex algebras and
their $\phi$-coordinated quasi modules. In Section 4, we study the
second deformed Heisenberg algebra in terms of quantum vertex
algebras and their $\phi$-coordinated quasi modules.

\section{Quantum vertex
algebras and their $\phi$-coordinated quasi modules}

This is a preliminary section. As we need, we recall the basic
notions and results, including the notions of (weak) quantum vertex
algebra and $\phi$-coordinated quasi module, and also including the
conceptual constructions.

We begin by recalling from \cite{li-qva1} the notion of weak quantum
vertex algebra, which generalizes the notion of vertex algebra and
that of vertex superalgebra.

\bd{dweak-qva} {\em A {\em weak quantum vertex algebra} is a vector
space $V$ equipped with a linear map
$$Y(\cdot,x): V\rightarrow \Hom(V,V((x)))\subset (\End
V)[[x,x^{-1}]]$$ and a vector ${\bf 1}\in V$, called the {\em vacuum
vector}, satisfying the following conditions: For $v\in V$,
$$Y({\bf 1},x)v=v,\ \ Y(v,x){\bf 1}\in V[[x]]\ \ \mbox{and }\
\lim_{x\rightarrow 0}Y(v,x){\bf 1}=v;$$  For $u,v,w\in V$, there
exists a nonnegative integer $l$ such that
\begin{eqnarray}
(x_{0}+x_{2})^{l}Y(u,x_{0}+x_{2})Y(v,x_{2})w=(x_{0}+x_{2})^{l}Y(Y(u,x_{0})v,x_{2})w;
\end{eqnarray}
and for $u,v\in V$, there exist
$$u^{(i)},v^{(i)}\in V,\ \ f_{i}(x)\in \C((x))\ \ \mbox{ for
}i=1,\dots,r$$ such that
\begin{eqnarray}
(x_{1}-x_{2})^{k}Y(u,x_{1})Y(v,x_{2})
=(x_{1}-x_{2})^{k}\sum_{i=1}^{r}f_{i}(x_{2}-x_{1})Y(v^{(i)},x_{2})Y(u^{(i)},x_{1})
\end{eqnarray}
for some nonnegative integer $k$.} \ed

A {\em rational quantum Yang-Baxter operator} on a vector space $U$
is a linear operator
$$\S(x):\ U\otimes U\rightarrow U\otimes U\otimes \C((x))$$
satisfying the quantum Yang-Baxter equation
$$\S^{12}(x)\S^{13}(x+z)\S^{23}(z)=\S^{23}(z)\S^{13}(x+z)\S^{12}(x).$$
It is said to be {\em unitary} if
$$\S(x)\S^{21}(-x)=1,$$
where $\S^{21}(x)=\sigma \S(x)\sigma$ with $\sigma$ denoting the
flip operator on $U\otimes U$.

\bd{dqva} {\em A {\em quantum vertex algebra} is a weak quantum
vertex algebra $V$ equipped with a unitary rational quantum
Yang-Baxter operator $\S(x)$ on $V$, satisfying
\begin{eqnarray}
&&\S(x)({\bf 1}\otimes v)={\bf 1}\otimes v\ \ \ \mbox{ for }v\in
V,\label{esvacuum}\\
&&[\D\otimes 1, \S(x)]=-\frac{d}{dx}\S(x),\label{d1s}\\
&&Y(u,x)v=e^{x\D}Y(-x)\S(-x)(v\otimes u)\ \ \mbox{ for }u,v\in V,\\
&&\S(x_{1})(Y(x_{2})\otimes 1)=(Y(x_{2})\otimes
1)\S^{23}(x_{1})\S^{13}(x_{1}+x_{2}).\label{sy1}
\end{eqnarray}
We denote a quantum vertex algebra by a pair $(V,\S)$.} \ed

Note that this very notion is a slight modification of the same
named notion in \cite{li-qva1} and \cite{li-qva2} with extra axioms
(\ref{esvacuum}) and (\ref{sy1}).

As we need, we recall the following important notion due to Etingof
and Kazhdan (see \cite{ek}):

\bd{dkl-nondeg} {\em A weak quantum vertex algebra $V$ is said to be
{\em non-degenerate} if for every positive integer $n$, the linear
map
$$Z_{n}: V^{\otimes n}\otimes \C((x_{1}))\cdots ((x_{n}))\rightarrow
V((x_{1}))\cdots ((x_{n})),$$ defined by
$$Z_{n}(v^{(1)}\otimes \cdots\otimes v^{(n)}\otimes f)
=fY(v^{(1)},x_{1})\cdots Y(v^{(n)},x_{n}){\bf 1}$$ for
$v^{(1)},\dots, v^{(n)}\in V,\; f\in \C((x_{1}))\cdots ((x_{n}))$,
is injective.} \ed

The following result can be found in \cite{li-qva1} (cf. \cite{ek}):

\bp{pnon-degenerate} Let $V$ be a weak quantum vertex algebra.
Assume that $V$ is non-degenerate. Then there exists a linear map
$\S(x): V\otimes V\rightarrow V\otimes V\otimes \C((x))$, which is
uniquely determined by
\begin{eqnarray*}
Y(u,x)v=e^{x\D}Y(-x)\S(-x)(v\otimes u) \ \ \ \mbox{for }u,v\in V.
\end{eqnarray*}
Furthermore, $(V,\S)$ carries the structure of a quantum vertex
algebra and the following relation holds
\begin{eqnarray}
[1\otimes \D, \S(x)]=\frac{d}{dx}\S(x).
\end{eqnarray}
 \ep

\br{rnotions} {\em Note that a quantum vertex algebra was defined as
a pair $(V,\S)$. In view of Proposition \ref{pnon-degenerate}, the
term ``a non-degenerate quantum vertex algebra'' (without reference
to a quantum Yang-Baxter operator) is unambiguous. Furthermore, if
$V$ is of countable dimension over $\C$ and if $V$ as a $V$-module
is irreducible, then by Corollary 3.10 of \cite{li-qva2}, $V$ is
non-degenerate. In view of this, the term ``irreducible quantum
vertex algebra'' is also unambiguous.} \er

We recall from \cite{li-qva1} a conceptual construction of weak
quantum vertex algebras and their modules. Let $W$ be a general
vector space. Set
\begin{eqnarray}
\E(W)=\Hom (W,W((x)))\subset (\End W)[[x,x^{-1}]].
\end{eqnarray}
The identity operator on $W$, denoted by $1_{W}$, is a typical
element of $\E(W)$.

\bd{dquasi-local-more} {\em A subset $U$ of $\E(W)$ is said to be
{\em $\S$-local} if for any $a(x),b(x)\in U$, there exist
$$u^{(i)}(x),v^{(i)}(x)\in U,\ \ f_{i}(x)\in \C((x))\ \ \mbox{ for
}i=1,\dots,r$$ such that
\begin{eqnarray}\label{ekab}
(x_{1}-x_{2})^{k}a(x_{1})b(x_{2})
=(x_{1}-x_{2})^{k}\sum_{i=1}^{r}f_{i}(x_{2}-x_{1})u^{(i)}(x_{2})v^{(i)}(x_{1})
\end{eqnarray}
for some nonnegative integer $k$, where by convention
$f_{i}(x_{2}-x_{1})$ is understood as an element of
$\C((x_{2}))((x_{1}))$.} \ed

Let $U$ be an $\S$-local subset of $\E(W)$ and let $a(x),b(x)\in U$.
Notice that (\ref{ekab}) implies
\begin{eqnarray}\label{ekab-truncation}
(x_{1}-x_{2})^{k}a(x_{1})b(x_{2})\in \Hom (W,W((x_{1},x_{2}))).
\end{eqnarray}
Define $a(x)_{n}b(x)\in \E(W)$ for $n\in \Z$ in terms of generating
function
$$Y_{\E}(a(x),z)b(x)=\sum_{n\in \Z}a(x)_{n}b(x) z^{-n-1}$$
by
\begin{eqnarray}
Y_{\E}(a(x),z)b(x)
=z^{-k}\left((x_{1}-x)^{k}a(x_{1})b(x)\right)|_{x_{1}=x+z},
\end{eqnarray}
where $k$ is any nonnegative integer such that
(\ref{ekab-truncation}) holds.  It was shown that
\begin{eqnarray*}
&&a(x)_{n}b(x)\nonumber\\
&=&\Res_{x_{1}}\left((x_{1}-x)^{n}a(x_{1})b(x)
-(-x+x_{1})^{n}\sum_{i=1}^{r}f_{i}(x-x_{1})u^{(i)}(x)v^{(i)}(x_{1})\right),
\end{eqnarray*}
where $u^{(i)}(x),v^{(i)}(x)\in U,\ \ f_{i}(x)\in \C((x))$
$(i=1,\dots,r)$ such that (\ref{ekab}) holds.

An $\S$-local subspace $K$ of $\E(W)$ is said to be {\em
$Y_{\E}$-closed} if
$$a(x)_{n}b(x)\in K\ \ \mbox{ for }a(x),b(x)\in K,\; n\in \Z.$$

The following result was obtained in \cite{li-qva1} (Theorem 5.8):

\bt{tweak-qva} Let $U$ be an $\S$-local subset of $\E(W)$. Then
there exists a $Y_{\E}$-closed $\S$-local subspace of $\E(W)$,
containing $1_{W}$ and $U$. Denote by $\<U\>$ the smallest such
$\S$-local subspace. Then $(\<U\>,Y_{\E},1_{W})$ carries the
structure of a weak quantum vertex algebra and $W$ is a faithful
$\<U\>$-module with $Y_{W}(\alpha(x),z)=\alpha(z)$ for $\alpha(x)\in
\<U\>$. \et

Throughout this section, we let
$$\phi=\phi(x,z)=xe^{z}\in \C[[x,z]].$$

The following notion was introduced in \cite{li-phi-module}:

\bd{dphi-quasimodule} {\em Let $V$ be a weak quantum vertex algebra.
A {\em $\phi$-coordinated quasi $V$-module} is a vector space $W$
equipped with a linear map
$$Y_{W}(\cdot,x): V\rightarrow \Hom (W,W((x)))\subset (\End
W)[[x,x^{-1}]],$$ satisfying the conditions that
$$Y_{W}({\bf 1},x)=1_{W}\ \ (\mbox{the identity operator on }W)$$
and that for any $u,v\in V$, there exists a nonzero polynomial
$p(x_{1},x_{2})$ such that
\begin{eqnarray}
p(x_{1},x_{2})Y_{W}(u,x_{1})Y_{W}(v,x_{2})\in
\Hom(W,W((x_{1},x_{2})))
\end{eqnarray}
and
\begin{eqnarray}
p(x_{2}e^{x_{0}},x_{2})Y_{W}(Y(u,x_{0})v,x_{2})
=\left(p(x_{1},x_{2})Y_{W}(u,x_{1})Y_{W}(v,x_{2})\right)|_{x_{1}=x_{2}e^{x_{0}}}.
\end{eqnarray}}
\ed

As we need, next we recall the conceptual construction from
\cite{li-phi-module}. Let $W$ be a general vector space.

\bd{dquasi-local-more} {\em A subset $U$ of $\E(W)$ is said to be
{\em quasi local} (see \cite{li-gamma}) if for any $a(x),b(x)\in U$,
there exists a nonzero polynomial $p(x_{1},x_{2})$ such that
\begin{eqnarray}
p(x_{1},x_{2})a(x_{1})b(x_{2})=p(x_{1},x_{2})b(x_{2})a(x_{1}).
\end{eqnarray}
A subset $U$ is said to be {\em quasi $S_{trig}$-local} if for any
$a(x),b(x)\in U$, there exist
$$u^{(i)}(x),v^{(i)}(x)\in U,\ \ f_{i}(x)\in \C(x)\ \ \mbox{ for
}i=1,\dots,r$$ such that
\begin{eqnarray}\label{epabpfuv}
p(x_{1},x_{2})a(x_{1})b(x_{2})
=p(x_{1},x_{2})\sum_{i=1}^{r}\iota_{x_{2},x_{1}}(f_{i}(x_{1}/x_{2}))u^{(i)}(x_{2})v^{(i)}(x_{1})
\end{eqnarray}
for some nonzero polynomial $p(x_{1},x_{2})$, where $\C(x)$ denotes
the field of rational functions and $\iota_{x_{2},x_{1}}$ is the
canonical field embedding of $\C(x_{1},x_{2})$ into
$\C((x_{2}))((x_{1}))$.} \ed

Let $U$ be a quasi $\S_{trig}$-local subset of $\E(W)$ and let
$a(x),b(x)\in U$. Notice that (\ref{epabpfuv}) implies
\begin{eqnarray}
p(x_{1},x_{2})a(x_{1})b(x_{2})\in \Hom (W,W((x_{1},x_{2}))).
\end{eqnarray}
Define $a(x)_{n}^{e}b(x)\in \E(W)$ for $n\in \Z$ in terms of
generating function
$$Y_{\E}^{e}(a(x),z)b(x)=\sum_{n\in \Z}a(x)_{n}^{e}b(x) z^{-n-1}$$
by
\begin{eqnarray}
Y_{\E}^{e}(a(x),z)b(x)
=\iota_{x,z}(1/p(xe^{z},x))\left(p(x_{1},x)a(x_{1})b(x)\right)|_{x_{1}=xe^{z}}.
\end{eqnarray}
A quasi $S_{trig}$-local subspace $K$ of $\E(W)$ is said to be {\em
$Y_{\E}^{e}$-closed} if
$$a(x)_{n}^{e}b(x)\in K\ \ \ \mbox{ for }a(x),b(x)\in K,\ n\in
\Z.$$

The following was proved in \cite{li-phi-module}:

\bp{pconvert} Let $V$ be a $Y_{\E}^{e}$-closed quasi compatible
subspace of $\E(W)$. Suppose
$$a(x),b(x),u_{i}(x),v_{i}(x)\in V,\; 0\ne p(x)\in \C[x],\; q_{i}(x)\in
\C(x)\ \ (i=1,\dots,r)$$ satisfy
\begin{eqnarray}\label{econdition}
p(x_{1}/x_{2})a(x_{1})b(x_{2})
=\sum_{i=1}^{r}p(x_{1}/x_{2})\iota_{x_{2},x_{1}}(q_{i}(x_{1}/x_{2}))
u_{i}(x_{2})v_{i}(x_{1}).
\end{eqnarray}
 Then
\begin{eqnarray}\label{epfsum01}
& &p(e^{x_{1}-x_{2}})Y_{\E}^{e}(a(x),x_{1})Y_{\E}^{e}(b(x),x_{2})
\nonumber\\
&=&p(e^{x_{1}-x_{2}})\sum_{i=1}^{r}\iota_{x_{2},x_{1}}(q_{i}(e^{x_{1}-x_{2}}))
Y_{\E}^{e}(u_{i}(x),x_{2})Y_{\E}^{e}(v_{i}(x),x_{1}).
\end{eqnarray}
Furthermore, we have
\begin{eqnarray}\label{epfsum}
& &(x_{1}-x_{2})^{k}Y_{\E}^{e}(a(x),x_{1})Y_{\E}^{e}(b(x),x_{2})
\nonumber\\
&=&(x_{1}-x_{2})^{k}\sum_{i=1}^{r}\iota_{x_{2},x_{1}}(q_{i}(e^{x_{1}-x_{2}}))
Y_{\E}^{e}(u_{i}(x),x_{2})Y_{\E}^{e}(v_{i}(x),x_{1}),
\end{eqnarray}
where $k$ is the multiplicity of zero of $p(x)$ at $x=1$. \ep

\bt{trecall1} Let $U$ be a quasi $S_{trig}$-local subset of $\E(W)$.
There exists a $Y_{\E}^{e}$-closed quasi $S_{trig}$-local subspace
containing $U$ and $1_{W}$. Denote by $\<U\>_{e}$ the smallest such
subspace. Then $(\<U\>_{e},Y_{\E}^{e},1_{W})$ carries the structure
of a weak quantum vertex algebra and $W$ is a $\phi$-coordinated
quasi $\<U\>_{e}$-module where $Y_{W}(\cdot,x)$ is given by
$Y_{W}(\alpha(x),z)=\alpha(z)$ for $\alpha(x)\in \<U\>_{e}$. \et

Combining Theorem \ref{trecall1} with Proposition \ref{pconvert} we
immediately have:

\bc{cvertex-algebra} Let $U$ be a quasi local subset of $\E(W)$.
Then the weak quantum vertex algebra $\<U\>_{e}$ obtained in Theorem
\ref{trecall1} is a vertex algebra. \ec

\bd{dcovariantmodule} {\em Let $V$ be a weak quantum vertex algebra,
$G$ an automorphism group of $V$, and $\chi: G\rightarrow
\C^{\times}$ a group homomorphism. We say a $\phi$-coordinated quasi
$V$-module $(W,Y_{W})$ is {\em $(G,\chi)$-covariant} if
\begin{eqnarray}\label{ecovaraince}
Y_{W}(gu,x)=Y_{W}(u,\chi(g)x)\ \ \ \mbox{ for }g\in G,\ u\in V.
\end{eqnarray}
In case that $G$ is a subgroup of $\C^{\times}$ with $\chi$ the
embedding, we simply drop $\chi$ from the notion. } \ed

The following is a useful technical result:

\bl{lsimplefacts} Let $V$, $G$, and $\chi$ be given as in Definition
\ref{dcovariantmodule}. Assume that $(W,Y_{W})$ is a
$\phi$-coordinated quasi $V$-module and that $U$ is a $G$-submodule
of $V$ which generates $V$ such that (\ref{ecovaraince}) holds for
$g\in G,\ u\in U$. Then $(W,Y_{W})$ is a $(G,\chi)$-covariant
$\phi$-coordinated quasi $V$-module. \el

\begin{proof} Suppose that $u,v\in V$ satisfy
$$Y_{W}(gu,x)=Y_{W}(u,\chi(g)x),\ \ \ \
Y_{W}(gv,x)=Y_{W}(v,\chi(g)x)\ \ \mbox{ for }g\in G.$$ Note that
$$Y_{W}(gY(u,z)v,x)=Y_{W}(Y(gu,z)gv,x).$$
There exists a nonzero homogeneous polynomial $p(x_{1},x_{2})$ such
that
$$p(xe^{z},x)Y_{W}(Y(gu,z)gv,x)=\left(p(x_{1},x)Y_{W}(gu,x_{1})Y_{W}(gv,x)\right)|_{x_{1}=xe^{z}}.$$
Then
\begin{eqnarray*}
&&\left(p(x_{1},x)Y_{W}(gu,x_{1})Y_{W}(gv,x)\right)|_{x_{1}=xe^{z}}\\
&=&\left(p(x_{1},x)Y_{W}(u,\chi(g)x_{1})Y_{W}(v,\chi(g)x)\right)|_{x_{1}=xe^{z}}\\
&=&\left(p(\chi(g)^{-1}x_{1},x)Y_{W}(u,x_{1})Y_{W}(v,\chi(g)x)\right)|_{x_{1}=\chi(g)xe^{z}}\\
&=&p(xe^{z},x)Y_{W}(Y(u,z)v,\chi(g)x).
\end{eqnarray*}
Consequently, we get
$$p(xe^{z},x)Y_{W}(Y(gu,z)gv,x)=p(xe^{z},x)Y_{W}(Y(u,z)v,\chi(g)x).$$
It was proved in \cite{li-phi-module} that $p(xe^{z},z)$ is a
nonzero element of $\C[[x,z]]\subset \C((x))((z))$. Multiplying both
sides by the inverse of $p(xe^{z},x)$ in $\C((x))[[z]]$ we obtain
$$Y_{W}(Y(gu,z)gv,x)=Y_{W}(Y(u,z)v,\chi(g)x).$$
Since $U$ generates $V$, it now follows that (\ref{ecovaraince})
holds for $g\in G,\ u\in V$.
\end{proof}

\section{$q$-Heisenberg Lie algebras}
In this section we study a certain $q$-Heisenberg Lie algebra in
terms of vertex algebras. We explicitly construct some vertex
algebras and then associate these vertex algebras and their
$\phi$-coordinated quasi modules to the Lie algebra.

We consider the following $q$-analog of the standard Heisenberg Lie
algebra.

\bd{dfirst} {\em Let $q$ be a nonzero complex number with $q\ne \pm
1$. Denote by $H_{q}$ the Lie algebra with generators $c$ and
$\beta_{m}$ $(m\in \Z)$, where $c$ is central, subject to relations
\begin{eqnarray}\label{efirst-relation}
[\beta_{m},\beta_{n}]=[m]_{q}\delta_{m+n,0}c,
\end{eqnarray}
where $[m]_{q}$ is the $q$-integer defined by
$$[m]_{q}=\frac{q^{m}-q^{-m}}{q-q^{-1}}.$$}
\ed

Form the generating function
\begin{eqnarray}
\beta(x)=\sum_{n\in \Z}\beta_{n}x^{-n}.
\end{eqnarray}
The defining relations (\ref{efirst-relation}) are equivalent to
\begin{eqnarray}\label{egood}
[\beta(x_{1}),\beta(x_{2})]=\frac{1}{q-q^{-1}}
\left(\delta\left(\frac{qx_{2}}{x_{1}}\right)-\delta\left(\frac{x_{2}}{qx_{1}}\right)\right)c.
\end{eqnarray}
 {}From this we get
\begin{eqnarray}
(x_{1}-qx_{2})(qx_{1}-x_{2})[\beta(x_{1}),\beta(x_{2})]=0.
\end{eqnarray}

\br{rdifferent-gf} {\em If we form the generating function
differently as $\beta(x)=\sum_{n\in \Z}\beta_{n}x^{-n-1}$, then we
have
\begin{eqnarray}
[\beta(x_{1}),\beta(x_{2})]=\frac{1}{q-q^{-1}}
\left(qx_{1}^{-2}\delta\left(\frac{qx_{2}}{x_{1}}\right)
-q^{-1}x_{1}^{-2}\delta\left(\frac{x_{2}}{qx_{1}}\right)\right)c,
\end{eqnarray}
which does not look as neat as (\ref{egood}).} \er

We say an $H_{q}$-module $W$ is of {\em level $\ell\in \C$} if $c$
acts on $W$ as scalar $\ell$, and we say an $H_{q}$-module $W$ is
{\em restricted} if for any $w\in W$, $\beta_{n}w=0$ for $n$
sufficiently large.

Next we associate this Lie algebra with vertex algebras in terms of
$\phi$-coordinated quasi modules. We first construct a vertex
algebra using another Lie algebra. Let $B$ be a vector space over
$\C$ with basis $\{ \beta^{(r)}\;|\;r\in \Z\}$. Equip $B$ with a
bilinear form defined by
\begin{eqnarray}\label{eform1}
\<\beta^{(r)},\beta^{(s)}\>=\frac{1}{q-q^{-1}}(\delta_{r,s+1}-\delta_{r,s-1})
=\frac{1}{q-q^{-1}}(\delta_{r,s+1}-\delta_{s,r+1})
\end{eqnarray}
for $r,s\in \Z$. It can be readily seen that this bilinear form is
skew-symmetric and non-degenerate. To the pair
$(B,\<\cdot,\cdot\>)$, we associate a Heisenberg Lie algebra
$$\hat{B}=B\otimes \C[t,t^{-1}]\oplus \C c,$$
where $c$ is central and
\begin{eqnarray}\label{eambn}
[a\otimes t^{m},b\otimes t^{n}]=\<a,b\>\delta_{m+n+1,0}c
\end{eqnarray}
for $a,b\in B,\; m,n\in \Z$.

For $b\in B$, form a generating function
\begin{eqnarray}
b(x)=\sum_{n\in \Z}b_{n}x^{-n-1},
\end{eqnarray}
where $b_{n}=b\otimes t^{n}$. Now, the defining relations
(\ref{eambn}) amount to
\begin{eqnarray}
[a(x_{1}),b(x_{2})]=\<a,b\>x_{1}^{-1}\delta\left(\frac{x_{2}}{x_{1}}\right)c
\end{eqnarray}
for $a,b\in B$.

Set
$$\hat{B}_{\ge 0}=(B\otimes \C[t])\oplus \C c,\ \ \ \
\hat{B}_{<0}=B\otimes t^{-1}\C[t^{-1}].$$ We see that $\hat{B}_{\ge
0}$ and $\hat{B}_{<0}$ are Lie subalgebras and $\hat{B}=\hat{B}_{\ge
0}\oplus \hat{B}_{<0}$ as a vector space. Let $\ell\in \C$ and
denote by $\C_{\ell}$ the one-dimensional $\hat{B}_{\ge 0}$-module
with $c$ acting as scalar $\ell$ and with $B\otimes \C[t]$ acting
trivially. Form the induced module
$$V_{\hat{B}}(\ell,0)=U(\hat{B})\otimes _{U(\hat{B}_{\ge
0})}\C_{\ell}.$$ Set ${\bf 1}=1\otimes 1\in V_{\hat{B}}(\ell,0)$.
Identify $B$ as a subspace of $V_{\hat{B}}(\ell,0)$ through the
linear map $b\in B\mapsto b(-1){\bf 1}$. It is now well known that
there exists a vertex algebra structure on $V_{\hat{B}}(\ell,0)$,
which is uniquely determined by the condition that ${\bf 1}$ is the
vacuum vector and $Y(b,x)=b(x)$ for $b\in B$. Furthermore, $B$ is a
generating subspace of vertex algebra $V_{\hat{B}}(\ell,0)$, where
for $r,s\in \Z$ and for $n\ge 0$,
\begin{eqnarray}
\beta^{(r)}_{n}\beta^{(s)}=\beta^{(r)}_{n}\beta^{(s)}_{-1}{\bf 1}
=\ell \<\beta^{(r)},\beta^{(s)}\>\delta_{n,0}{\bf 1}
\end{eqnarray}
as $\beta^{(r)}_{n}{\bf 1}=0$. It is clear that
$V_{\hat{B}}(\ell,0)$ is an irreducible $\hat{B}$-module. It follows
that vertex algebra $V_{\hat{B}}(\ell,0)$ is simple.

Furthermore we have:

\bl{lZ-action} The abelian group $\Z$ acts on $V_{\hat{B}}(\ell,0)$
as an automorphism group such that
$$\rho_{m}(\beta^{(r)})=\beta^{(r+m)}\ \ \ \mbox{ for }m,r\in \Z.$$
\el

\begin{proof} For $m\in \Z$, let $\hat{\rho}_{m}$ be the linear
automorphism of $B$ defined by
$\hat{\rho}_{m}(\beta^{(n)})=\beta^{(m+n)}$ for $n\in \Z$. It is
clear that $\hat{\rho}_{m}$ preserves the bilinear form
$\<\cdot,\cdot\>$ (see (\ref{eform1})). It follows that
$\hat{\rho}_{m}$ gives an automorphism of the Lie algebra $H_{q}$
such that $\hat{\rho}_{m} (c)=c$ and
$$\hat{\rho}_{m} (\beta^{(n)}\otimes t^{r})=\beta^{(m+n)}\otimes
t^{r}\ \ \ \mbox{ for }n,r\in \Z.$$ This gives rise to an
automorphism of $U(H_{q})$, also denoted by $\hat{\rho}_{m}$.
Furthermore, $\hat{\rho}_{m}$ gives rise to a linear automorphism
$\rho_{m}$ of $V_{\hat{B}}(\ell,0)$ such that $\rho_{m}({\bf
1})={\bf 1}$ and
$$\rho_{m}(av)=\hat{\rho}_{m}(a)\rho_{m}(v)\ \ \ \mbox{ for }a\in
U(H_{q}),\; v\in V_{\hat{B}}(\ell,0).$$ As $V_{\hat{B}}(\ell,0)$ as
a vertex algebra is generated by $B$, it follows that $\rho_{m}$ is
an automorphism of $V_{\hat{B}}(\ell,0)$ viewed as a vertex algebra.
This gives an action of $\Z$ on $V_{\hat{B}}(\ell,0)$ by
automorphisms.
\end{proof}

Set
$$\Gamma_{q}=\{ q^{n}\;|\; n\in \Z\} \subset \C^{\times}.$$
Recall that a $\phi$-coordinated quasi $V_{\hat{B}}(\ell,0)$-module
$(W,Y_{W})$ is $\Gamma_{q}$-covariant if
\begin{eqnarray}
Y_{W}(\rho_{m}(v),x)=Y_{W}(v,q^{m}x) \ \ \ \mbox{ for }v\in
V_{\hat{B}}(\ell,0),\; m\in \Z.
\end{eqnarray}

As our main result of this section we have:

\bt{tmail-1} Let $W$ be a restricted $H_{q}$-module of level $\ell$.
Then there exists a $\Gamma_{q}$-covariant $\phi$-coordinated quasi
$V_{\hat{B}}(\ell,0)$-module structure on $W$, which is uniquely
determined by $Y_{W}(\beta^{(r)},x)=\beta(q^{r}x)$ for $r\in \Z$.
\et

\begin{proof} Set $U_{W}=\{ \beta(q^{r}x)\;|\; r\in \Z\}\subset \E(W)$.
For $r,s\in \Z$, we have
\begin{eqnarray}\label{ebqrbqs}
[\beta(q^{r}x_{1}),\beta(q^{s}x_{2})]=\frac{1}{q-q^{-1}}
\left(\delta\left(q^{s-r+1}\frac{x_{2}}{x_{1}}\right)
-\delta\left(q^{s-r-1}\frac{x_{2}}{x_{1}}\right)\right)\ell,
\end{eqnarray}
which implies
\begin{eqnarray}\label{ebqrbqs=0}
\left(\frac{x_{1}}{x_{2}}-q^{s-r+1}\right)\left(\frac{x_{1}}{x_{2}}-q^{s-r-1}\right)
[\beta(q^{r}x_{1}),\beta(q^{s}x_{2})]=0.
\end{eqnarray}
Then $U_{W}$ is a quasi local subset of $\E(W)$. By Corollary
\ref{cvertex-algebra}, $U_{W}$ generates a vertex algebra $V_{W}$
and $W$ is a $\phi$-coordinated quasi $V_{W}$-module where
$Y_{W}(a(x),z)=a(z)$ for $a(x)\in V_{W}$.

With (\ref{ebqrbqs=0}),  by Lemma 6.7 of \cite{li-phi-module} we
have
\begin{eqnarray*}
&&\beta(q^{r}x)_{n}^{e}\beta(q^{s}x)=0 \ \ \  \mbox{ for
}r-s\ne \pm 1,\; n\ge 0,\nonumber\\
 &&\beta(q^{r}x)_{n}^{e}\beta(q^{s}x)=0\ \ \ \mbox{ for }r-s=\pm 1,\; n\ge 1.
\end{eqnarray*}
If $r-s=1$, also by Lemma 6.7 of \cite{li-phi-module} we have
\begin{eqnarray*}
&&(1-q^{-2})\beta(q^{r}x_{1})_{0}^{e}\beta(q^{s}x_{2})\\
&=&\Res_{x_{1}}\frac{1}{x_{1}-x}\left(\frac{x_{1}}{x}-1\right)
\left(\frac{x_{1}}{x}-q^{-2}\right)\beta(q^{r}x_{1})\beta(q^{s}x)\\
&&-\Res_{x_{1}}\frac{1}{-x+x_{1}}\left(\frac{x_{1}}{x}-1\right)
\left(\frac{x_{1}}{x}-q^{-2}\right)\beta(q^{s}x)\beta(q^{r}x_{1})\\
&=&\Res_{x_{1}}x^{-2}(x_{1}-q^{-2}x)[\beta(q^{r}x_{1}),\beta(q^{s}x)]\\
&=&\frac{1}{q-q^{-1}}\Res_{x_{1}}x^{-2}(x_{1}-q^{-2}x)
\left(\delta\left(\frac{x}{x_{1}}\right)
-\delta\left(q^{-2}\frac{x}{x_{1}}\right)\right)\ell\\
 &=&\ell \frac{1-q^{-2}}{q-q^{-1}}.
\end{eqnarray*}
That is,
$$\beta(q^{r}x_{1})_{0}^{e}\beta(q^{s}x_{2})=\frac{\ell}{q-q^{-1}}1_{W}.$$
Similarly, if $r-s=-1$, we have
$$\beta(q^{r}x_{1})_{0}^{e}\beta(q^{s}x_{2})=-\frac{\ell}{q-q^{-1}}1_{W}.$$
Consequently, we have
\begin{eqnarray}
[Y_{\E}^{\phi}(\beta(q^{r}x),x_{1}),Y_{\E}^{\phi}(\beta(q^{s}x),x_{2})]
=\frac{\ell}{q-q^{-1}}
(\delta_{r,s+1}-\delta_{r,s-1})x_{1}^{-1}\delta\left(\frac{x_{2}}{x_{1}}\right).
\end{eqnarray}
Thus $V_{W}$ is a restricted $\hat{B}$-module of level $\ell$ with
$\beta^{(r)}_{n}=\beta(q^{r}x)_{n}$ for $r,n\in \Z$. Furthermore,
$V_{W}$ together with vector $1_{W}$ is a vacuum $\hat{B}$-module of
level $\ell$. It follows from the construction of
$V_{\hat{B}}(\ell,0)$ that there exists a $\hat{B}$-module
homomorphism $\psi$ from $V_{\hat{B}}(\ell,0)$ to $V_{W}$ with
$\psi({\bf 1})=1_{W}$. For $r\in \Z$, we have
$$\psi(\beta^{(r)})=\psi(\beta^{(r)}_{-1}{\bf 1})
=\beta(q^{r}x)_{-1}1_{W}=\beta(q^{r}x)\in V_{W}.$$ It follows that
$\psi$ is a homomorphism of vertex algebras. Consequently, $W$ is a
$\phi$-coordinated quasi $V_{\hat{B}}(\ell,0)$-module. For $m,r\in
\Z$, we have
$$Y_{W}(\rho_{m}\beta^{(r)},z)=Y_{W}(\beta^{(m+r)},z)=\beta(q^{m+r}z)
=Y_{W}(\beta^{(r)},q^{m}z).$$  The covariance property follows from
Lemma \ref{lsimplefacts}.
\end{proof}

\section{Heisenberg Lie algebra $\widetilde{{\mathcal{H}}}_{q}$}
In this section, we study another deformed Heisenberg Lie algebra.
We explicitly construct quantum vertex algebras and we then
associate these quantum vertex algebras and their $\phi$-coordinated
quasi modules to the deformed Heisenberg Lie algebra.

\bd{dhq} {\em Let $q$ be a complex number. Denote by
$\widetilde{{\mathcal{H}}}_{q}$ the Lie algebra with generators
$\beta_{n}$ $(n\in \Z)$ and $c$, with $c$ central, subject to
relations
\begin{eqnarray}\label{esecond-relation}
[\beta_{m},\beta_{n}]=m(1-q^{|m|})\delta_{m+n,0}c
\end{eqnarray}
for $m,n\in \Z$.} \ed

If $q=\pm 1$, $\widetilde{{\mathcal{H}}}_{q}$ is just an abelian Lie
algebra. If $q=0$, this becomes the standard free field Heisenberg
Lie algebra. For the rest of this section, we assume that {\em $q$
is neither zero nor a root of unity.}

Form the generating function
\begin{eqnarray}
\beta(x)=\sum_{n\in \Z}\beta_{n}x^{-n}.
\end{eqnarray}

\bl{lrelation-generating} The defining relations
(\ref{esecond-relation}) amount to
\begin{eqnarray}\label{erelation-generating}
[\beta(x_{1}),\beta(x_{2})]=\left(\left(x_{2}\frac{\partial}{\partial
x_{2}}\right)\delta\left(\frac{x_{2}}{x_{1}}\right)
+\frac{qx_{1}/x_{2}}{(1-qx_{1}/x_{2})^{2}}
-\frac{qx_{2}/x_{1}}{(1-qx_{2}/x_{1})^{2}}\right)c.
\end{eqnarray}
Furthermore, the latter is equivalent to
\begin{eqnarray}\label{eqlocality-commutator}
&&\beta(x_{1})\beta(x_{2})-\beta(x_{2})\beta(x_{1})=
\frac{qx_{1}/x_{2}}{(1-qx_{1}/x_{2})^{2}}c
-\frac{q^{-1}x_{1}/x_{2}}{(1-q^{-1}x_{1}/x_{2})^{2}}c\hspace{1cm}\nonumber\\
&&\hspace{2cm}+\left(x_{2}\frac{\partial}{\partial
x_{2}}\right)\left(
\delta\left(\frac{x_{2}}{x_{1}}\right)-\delta\left(\frac{qx_{2}}{x_{1}}\right)\right)c.
\end{eqnarray}
\el

\begin{proof} It is straightforward:
\begin{eqnarray*}
&&[\beta(x_{1}),\beta(x_{2})]=\sum_{m\in
\Z}m(1-q^{|m|})x_{1}^{-m}x_{2}^{m}c\\
&=&\sum_{m\in \Z}mx_{1}^{-m}x_{2}^{m}c-\sum_{m\ge
0}mq^{m}x_{1}^{-m}x_{2}^{m}c-\sum_{m<0}mq^{-m}x_{1}^{-m}x_{2}^{m}c\\
&=&\left(x_{2}\frac{\partial}{\partial x_{2}}\right)
\left(\sum_{m\in \Z}x_{1}^{-m}x_{2}^{m}-\sum_{m\ge
0}q^{m}x_{1}^{-m}x_{2}^{m}-\sum_{m\le 0}q^{-m}x_{1}^{-m}x_{2}^{m}\right)c\\
&=&\left(x_{2}\frac{\partial}{\partial x_{2}}\right)
\left(\sum_{m\in \Z}x_{1}^{-m}x_{2}^{m}-\sum_{m\ge 0}q^{m}x_{1}^{-m}x_{2}^{m}
-\sum_{m\ge 0}q^{m}x_{1}^{m}x_{2}^{-m}\right)c\\
&=&\left(x_{2}\frac{\partial}{\partial
x_{2}}\right)\left(\delta\left(\frac{x_{2}}{x_{1}}\right)
-\frac{1}{1-qx_{2}/x_{1}}-\frac{1}{1-qx_{1}/x_{2}}\right)c\\
&=&\left(\left(x_{2}\frac{\partial}{\partial
x_{2}}\right)\delta\left(\frac{x_{2}}{x_{1}}\right)
-\frac{qx_{2}/x_{1}}{(1-qx_{2}/x_{1})^{2}}+\frac{qx_{1}/x_{2}}{(1-qx_{1}/x_{2})^{2}}\right)c.
\end{eqnarray*}
Using the calculation above we obtain
\begin{eqnarray*}
&&[\beta(x_{1}),\beta(x_{2})]\\
&=&\left(x_{2}\frac{\partial}{\partial x_{2}}\right)
\left(\sum_{m\in \Z}x_{1}^{-m}x_{2}^{m}-\sum_{m\ge
0}q^{m}x_{1}^{-m}x_{2}^{m}-\sum_{m\le 0}q^{-m}x_{1}^{-m}x_{2}^{m}\right)c\\
&=&\left(x_{2}\frac{\partial}{\partial
x_{2}}\right)\left(\delta\left(\frac{x_{2}}{x_{1}}\right)
-\delta\left(\frac{qx_{2}}{x_{1}}\right)
+\sum_{m<0}q^{m}x_{1}^{-m}x_{2}^{m}-\sum_{m\le 0}q^{-m}x_{1}^{-m}x_{2}^{m}\right)c\\
&=&\left(x_{2}\frac{\partial}{\partial
x_{2}}\right)\left(\delta\left(\frac{x_{2}}{x_{1}}\right)
-\delta\left(\frac{qx_{2}}{x_{1}}\right)+\frac{x_{1}/(qx_{2})}{1-x_{1}/(qx_{2})}
-\frac{1}{1-qx_{1}/x_{2}}\right)c\\
&=&\left(x_{2}\frac{\partial}{\partial
x_{2}}\right)\left(\delta\left(\frac{x_{2}}{x_{1}}\right)
-\delta\left(\frac{qx_{2}}{x_{1}}\right)\right)c-
\frac{q^{-1}x_{1}/x_{2}}{(1-q^{-1}x_{1}/x_{2})^{2}}c
+\frac{qx_{1}/x_{2}}{(1-qx_{1}/x_{2})^{2}}c,
\end{eqnarray*}
as desired.
\end{proof}

{}From (\ref{eqlocality-commutator}) we have
\begin{eqnarray}\label{eqlocality}
&&(x_{1}-x_{2})^{2}(x_{1}-qx_{2})^{2}\beta(x_{1})\beta(x_{2})
\nonumber\\
&=&(x_{1}-x_{2})^{2}(x_{1}-qx_{2})^{2}
\left(\beta(x_{2})\beta(x_{1})+
\frac{qx_{1}/x_{2}}{(1-qx_{1}/x_{2})^{2}}c
-\frac{q^{-1}x_{1}/x_{2}}{(1-q^{-1}x_{1}/x_{2})^{2}}c \right).\ \ \
\ \ \
\end{eqnarray}

\bd{drestricted} {\em We say an $\widetilde{\mathcal{H}}_{q}$-module
$W$ is {\em restricted} if $\beta(x)\in \E(W)$, i.e., for any $w\in
W$, $\beta_{n}w=0$ for $n$ sufficiently large. If $c$ acts on $W$ as
a scalar $\ell\in \C$ we say $W$ is of {\em level} $\ell$.} \ed

Let $W$ be a restricted $\widetilde{\mathcal{H}}_{q}$-module of
level $\ell$. Set $$U_{W}=\{ \beta(x),1_{W}\}\subset \E(W).$$ {}From
(\ref{eqlocality}) we see that $U_{W}$ is a quasi $\S_{trig}$-local
subset of $\E(W)$. By Theorem \ref{trecall1}, $U_{W}$ generates a
weak quantum vertex algebra $\<U_{W}\>_{e}$ with $W$ as a
$\phi$-coordinated quasi module. Next, we determine this weak
quantum vertex algebra.

\bd{dbhat} {\em Let $\hat{B}_{q}$ denote the Lie algebra with
generators $\beta^{(r)}_{m}$ $(r,m\in \Z)$ and $c$, with $c$
central, subject to relations
\begin{eqnarray}
&&\beta^{(r)}(x_{1})\beta^{(s)}(x_{2})
-\beta^{(s)}(x_{2})\beta^{(r)}(x_{1})\nonumber\\
&=&
\iota_{x_{2},x_{1}}\left(\frac{q^{r-s+1}e^{x_{1}-x_{2}}}{(1-q^{r-s+1}e^{x_{1}-x_{2}})^{2}}-
\frac{q^{r-s-1}e^{x_{1}-x_{2}}}{(1-q^{r-s-1}e^{x_{1}-x_{2}})^{2}}\right)c\nonumber\\
&&+(\delta_{r,s}-\delta_{r,s+1})\frac{\partial}{\partial x_{2}}
x_{1}^{-1}\delta\left(\frac{x_{2}}{x_{1}}\right)c
\end{eqnarray}
for $r,s\in \Z$, where
\begin{eqnarray}
\beta^{(r)}(x)=\sum_{m\in \Z}\beta^{(r)}_{m}x^{-m-1}.
\end{eqnarray}} \ed

\br{rffunction} {\em For $n\in \Z$, set
\begin{eqnarray}
f_{n}(x)=\iota_{x,0}\left(\frac{q^{n+1}e^{x}}{(1-q^{n+1}e^{x})^{2}}-
\frac{q^{n-1}e^{x}}{(1-q^{n-1}e^{x})^{2}}\right)\in \C((x)),
\end{eqnarray}
where $\iota_{x,0}(\cdot)$ stands for the formal Laurent series
expansion at $x=0$. Then
$$\iota_{x_{2},x_{1}}\left(\frac{q^{r-s+1}e^{x_{1}-x_{2}}}{(1-q^{r-s+1}e^{x_{1}-x_{2}})^{2}}-
\frac{q^{r-s-1}e^{x_{1}-x_{2}}}{(1-q^{r-s-1}e^{x_{1}-x_{2}})^{2}}\right)
=f_{r-s}(-x_{2}+x_{1})$$ for $r,s\in \Z$. We see that
\begin{eqnarray}\label{efnx-exp}
f_{n}(x)=\delta_{n+1,0}x^{-2}-\delta_{n-1,0}x^{-2}+O(1).
\end{eqnarray}} \er

Let $\hat{B}_{q}^{+}$ denote the subspace spanned by
$\beta^{(r)}_{n}$ for $r\in \Z,\ n\ge 0$. Noticing that
$$\iota_{x_{2},x_{1}}\left(\frac{q^{r-s+1}e^{x_{1}-x_{2}}}{(1-q^{r-s+1}e^{x_{1}-x_{2}})^{2}}-
\frac{q^{r-s-1}e^{x_{1}-x_{2}}}{(1-q^{r-s-1}e^{x_{1}-x_{2}})^{2}}\right)=f_{r-s}(-x_{2}+x_{1})$$
involves only nonnegative powers of $x_{1}$, we get
\begin{eqnarray}\label{ermsnhalf}
[\beta^{(r)}_{m},\beta^{(s)}_{n}]=(\delta_{r,s}-\delta_{r,s+1})m\delta_{m+n,0}c
\end{eqnarray}
for $r,s,m,n\in \Z$ with $m\ge 0$. It follows that $\hat{B}_{q}^{+}$
is an abelian subalgebra.

A $\hat{B}_{q}$-module $W$ is said to be {\em restricted} if
$\beta^{(r)}(x)\in \E(W)$ for $r\in \Z$, i.e., for any $w\in W,\
r\in \Z$, $\beta^{(r)}_{n}w=0$ for $n$ sufficiently large. A
$\hat{B}_{q}$-module $W$ is said to be of {\em level} $\ell\in \C$
if $c$ acts on $W$ as scalar $\ell$. A {\em vacuum}
$\hat{B}_{q}$-module of level $\ell$ is a module $W$ equipped with a
vector $w_{0}\in W$ such that $W$ is cyclic on $w_{0}$ and
$$\beta^{(r)}_{n}w_{0}=0\ \ \ \mbox{ for }r,n\in \Z\ \mbox{with
}n\ge 0.$$

\br{rvacuum-restricted} {\em Let $(W,w_{0})$ be a vacuum
$\hat{B}_{q}$-module. By definition we have $W=U(\hat{B}_{q})w_{0}$.
On the other hand, it follows from (\ref{ermsnhalf}) (and induction)
that for any $r\in \Z$ and for any $u\in U(\hat{B}_{q})$,
$\beta^{(r)}_{m}u=u\beta^{(r)}_{m}$ for $m$ sufficiently large. Then
it follows that $W$ is restricted.} \er

Let $\ell\in \C$. Denote by $\C_{\ell}=\C$ the $1$-dimensional
$(\hat{B}_{q}^{+}+\C c)$-module with $\hat{B}_{q}^{+}$ acting
trivially and with $c$ acting as scalar $\ell$. Form the induced
module
\begin{eqnarray}
V_{\hat{B}_{q}}(\ell,0)=U(\hat{B}_{q})\otimes_{U(\hat{B}_{q}^{+}+\C
c)} \C_{\ell}.
\end{eqnarray}
Set ${\bf 1}=1\otimes 1\in V_{\hat{B}_{q}}(\ell,0)$. Identify
$\beta^{(r)}$ with $\beta^{(r)}_{-1}{\bf 1}$ for $r\in \Z$. It is
clear that $V_{\hat{B}_{q}}(\ell,0)$ is a vacuum
$\hat{B}_{q}$-module of level $\ell$, which is universal in the
obvious sense.

\bt{tmain-2} Let $\ell\in\C$. There exists a weak quantum vertex
algebra structure on $V_{\hat{B}_{q}}(\ell,0)$, which is uniquely
determined by the condition that ${\bf 1}$ is the vacuum vector and
$Y(\beta^{(r)},x)=\beta^{(r)}(x)$ for $r\in \Z$. Furthermore, if
$\ell\ne 0$, $V_{\hat{B}_{q}}(\ell,0)$ is an irreducible quantum
vertex algebra. On the other hand, for any restricted
$\hat{B}_{q}$-module $W$ of level $\ell$, there exists a
$V_{\hat{B}_{q}}(\ell,0)$-module structure uniquely determined by
the condition that $Y_{W}(\beta^{(r)},x)=\beta^{(r)}(x)$ for $r\in
\Z$. \et

\begin{proof} Let $W$ be any restricted $\hat{B}_{q}$-module of
level $\ell$. Set $$U_{W}=\{1_{W}\}\cup \{ \beta^{(r)}(x)\;|\; r\in
\Z\}.$$ {}From the defining relations of $\hat{B}_{q}$, we see that
$U_{W}$ is an $\S$-local subset of $\E(W)$. By Theorem 2.7, $U_{W}$
generates a weak quantum vertex algebra $\<U_{W}\>$. Furthermore, it
follows from Proposition 6.6 of \cite{li-qva1} that
\begin{eqnarray}
&&Y_{\E}(\beta^{(r)}(x),x_{1})Y_{\E}(\beta^{(s)}(x),x_{2})
-Y_{\E}(\beta^{(s)}(x),x_{2})Y_{\E}(\beta^{(r)}(x),x_{1})\nonumber\\
&=&
\iota_{x_{2},x_{1}}\left(\frac{q^{r-s+1}e^{x_{1}-x_{2}}}{(1-q^{r-s+1}e^{x_{1}-x_{2}})^{2}}-
\frac{q^{r-s-1}e^{x_{1}-x_{2}}}{(1-q^{r-s-1}e^{x_{1}-x_{2}})^{2}}\right)\ell 1_{W}\nonumber\\
&&\hspace{1cm}+\ell(\delta_{r,s}-\delta_{r,s+1})\frac{\partial}{\partial
x_{2}} x_{1}^{-1}\delta\left(\frac{x_{2}}{x_{1}}\right)1_{W}
\end{eqnarray}
for $r,s\in \Z$. {}From this we see that $\<U_{W}\>$ is a vacuum
$\hat{B}_{q}$-module of level $\ell$ with $\beta^{(r)}(z)$ acting as
$Y_{\E}(\beta^{(r)}(x),z)$ for $r\in \Z$. Then it follows from the
universality of $V_{\hat{B}_{q}}(\ell,0)$ that there exists a
$\hat{B}_{q}$-module homomorphism $\psi$ from
$V_{\hat{B}_{q}}(\ell,0)$ to $\<U_{W}\>$, sending ${\bf 1}$ to
$1_{W}$. For $r\in \Z$, we have
$$\psi(\beta^{(r)}_{n}v)=\beta^{(r)}(x)_{n}\psi(v)\ \ \mbox{ for }v\in
V_{\hat{B}_{q}}(\ell,0),$$ in particular,
$$\psi(\beta^{(r)})=\psi(\beta^{(r)}_{-1}{\bf
1})=\beta^{(r)}(x)_{-1}1_{W}=\beta^{(r)}(x).$$ Specializing $W$ to
$V_{\hat{B}_{q}}(\ell,0)$, by Theorem 2.9 of \cite{li-qva2} with
$U=\{\beta^{(r)}\;|\; r\in \Z\}$ and
$$Y_{0}(\beta^{(r)},x)=\beta^{(r)}(x)\ \ \ \mbox{ for }r\in \Z,$$
we have a weak quantum vertex algebra structure on
$V_{\hat{B}_{q}}(\ell,0)$ with all the requirements and such a
structure is unique.

Come back to a general restricted $\hat{B}_{q}$-module $W$ of level
$\ell$. We have a $\hat{B}_{q}$-module homomorphism $\psi$ from
$V_{\hat{B}_{q}}(\ell,0)$ to $\<U_{W}\>$, sending ${\bf 1}$ to
$1_{W}$. For $r\in \Z$, we have
$$\psi(Y(\beta^{(r)},z)v)=\psi(\beta^{(r)}(z)v)=Y_{\E}(\beta^{(r)}(x),z)\psi(v)
=Y_{\E}(\psi(\beta^{(r)}),z)\psi(v)$$ for $v\in
V_{\hat{B}_{q}}(\ell,0)$. Since $\{\beta^{(r)}\;|\; r\in \Z\}$
generates $V_{\hat{B}_{q}}(\ell,0)$, it follows that $\psi$ is a
homomorphism of weak quantum vertex algebras. As $W$ is a canonical
$\<U_{W}\>$-module, $W$ becomes a $V_{\hat{B}_{q}}(\ell,0)$-module
 through the homomorphism $\psi$, where
$$Y_{W}(\beta^{(r)},z)=Y_{W}(\psi(\beta^{(r)}),z)=Y_{W}(\beta^{(r)}(x),z)=\beta^{(r)}(z)$$
for $r\in \Z$. The uniqueness is clear.

Next, we prove that $V_{\hat{B}_{q}}(\ell,0)$ is an irreducible
$\hat{B}_{q}$-module. First we define a descending filtration for
Lie algebra $\hat{B}_{q}$. For $n\ge 1$, set
$$\hat{B}_{q}[n]=\mbox{span}\{ \beta^{(r)}_{m}\;|\; r,m\in \Z,\ m\ge n\},$$
and for $n\le 0$, set
$$\hat{B}_{q}[n]=\mbox{span}\{\beta^{(r)}_{m}\;|\; r,m\in \Z,\ m\ge n\}+\C c.$$
Using (\ref{ermsnhalf}) we get
$$\left[\hat{B}_{q}[m],\hat{B}_{q}[n]\right]\subset \hat{B}_{q}[m+n]
\ \ \ \mbox{ for }m,n\in \Z.$$ Denote the associated graded Lie
algebra by $L$. For $r,m\in \Z$, set
$$\bar{\beta}^{(r)}_{m}=\beta^{(r)}_{m}+\hat{B}_{q}[m+1]\in L$$
and set
$$\bar{c}=c+\hat{B}_{q}[1]\in L.$$
With (\ref{efnx-exp}) we see that $L$ satisfies relations
\begin{eqnarray}
&&[\bar{\beta}^{(r)}(x_{1}),\bar{\beta}^{(s)}(x_{2})]\nonumber\\
&=&(\delta_{r,s}-\delta_{r,s+1})\frac{\partial}{\partial x_{2}}
x_{1}^{-1}\delta\left(\frac{x_{2}}{x_{1}}\right)\bar{c}
+\left(\delta_{r-s+1,0}-\delta_{r-s-1,0}\right)(x_{2}-x_{1})^{-2}\bar{c}\nonumber\\
 &=&\delta_{r,s}\frac{\partial}{\partial x_{2}}
x_{1}^{-1}\delta\left(\frac{x_{2}}{x_{1}}\right)\bar{c}
+\delta_{r,s-1}\frac{1}{(x_{2}-x_{1})^{2}}\bar{c}
-\delta_{r,s+1}\frac{1}{(x_{1}-x_{2})^{2}}\bar{c}
\end{eqnarray}
for $r,s\in \Z$. In terms of components we have
\begin{eqnarray}
[\bar{\beta}^{(r)}_{m},\bar{\beta}^{(s)}_{n}]
=\left(m\delta_{r,s}+\delta_{r,s-1}\frac{|m|-m}{2}
-\delta_{r,s+1}\frac{|m|+m}{2}\right)\delta_{m+n,0}\bar{c}
\end{eqnarray}
for $m,n\in \Z$. In Lemma \ref{linproof} below we shall prove that
for nonzero $\ell\in \C$, every (nonzero) vacuum $L$-module of level
$\ell$ is irreducible. Then by Proposition 2.11 of \cite{kl},
$V_{\hat{B}_{q}}(\ell,0)$ is an irreducible $\hat{B}_{q}$-module. It
follows that $V_{\hat{B}_{q}}(\ell,0)$ as a (left)
$V_{\hat{B}_{q}}(\ell,0)$-module is irreducible. Now, the proof is
complete.
\end{proof}

The following is the result we need at the end of the proof of
Theorem \ref{tmain-2}:

\bl{linproof} Let $L$ be a Heisenberg algebra with basis $\{c\}\cup
\{ \bar{\beta}^{(r)}_{m}\;|\; r,m\in \Z\}$, where $c$ is central and
\begin{eqnarray*}
[\bar{\beta}^{(r)}_{m},\bar{\beta}^{(s)}_{n}]
=\left(m\delta_{r,s}+\delta_{r,s-1}\frac{|m|-m}{2}
-\delta_{r,s+1}\frac{|m|+m}{2}\right)\delta_{m+n,0}c
\end{eqnarray*}
for $m,n,r,s\in \Z$. Suppose that $W$ is an $L$-module of nonzero
level $\ell$ equipped with a nonzero vector $w_{0}\in V$, satisfying
$$W=U(L)w_{0}\ \ \mbox{ and }\ \ \bar{\beta}^{(r)}_{m}w_{0}=0
\  \mbox{ for }r,m\in \Z \mbox{ with } m\ge 0.$$ Then $W$ is
irreducible.\el

\begin{proof} Set
$$V=\C[x_{n}^{(r)}\;|\; r\in \Z,\; n\ge 1],$$
a commutative polynomial algebra in variables $x_{n}^{(r)}$. For $r,
n\in \Z$ with $n\ge 1$, set
\begin{eqnarray*}
\bar{\beta}^{(r)}_{0}=0,\ \ \ \ \bar{\beta}^{(r)}_{-n}=x_{n}^{(r)},\
\ \ \ \bar{\beta}^{(r)}_{n}=\ell n\left(\frac{\partial}{\partial
x_{n}^{(r)}}-\frac{\partial}{\partial x_{n}^{(r-1)}}\right).
\end{eqnarray*}
 It is straightforward to check
that this defines an $L$-module structure on $V$ of level $\ell$.
Now we show that $V$ is an irreducible $L$-module. Let $A$ be a
nonzero $L$-submodule of $V$. Set $\deg x_{n}^{(r)}=1$ for $n,r\in
\Z$ with $n\ge 1$. Let $P\in A$ be a nonzero polynomial with least
degree. If $\deg P=0$, then $A=V$ as $V$ is clearly cyclic on $1$.
We next show that this must be the case. Assume $\deg P\ge 1$. There
exists $r\in \Z$ such that
$$\frac{\partial P}{\partial
x_{m}^{(s)}}=0\ \ \ \mbox{ for all }m\ge 1,\; s<r,$$ and
$\frac{\partial P}{\partial x_{n}^{(r)}}\ne 0$ for some $n\ge 1$.
Then
$$\bar{\beta}^{(r)}_{n}\cdot P=\ell n\left(\frac{\partial}{\partial
x_{n}^{(r)}}-\frac{\partial}{\partial x_{n}^{(r-1)}}\right)P=\ell
n\frac{\partial P}{\partial x_{n}^{(r)}}.$$ We have that
$\frac{\partial P}{\partial x_{n}^{(r)}}\in A$ and $\frac{\partial
P}{\partial x_{n}^{(r)}}\ne 0$, a contradiction. Therefore, $V$ is
irreducible. It follows from the standard argument (cf. \cite{ll})
that $W$ is isomorphic to $V$.
\end{proof}

The following is analogous to Lemma \ref{lZ-action}:

\bl{lautomorphisms} For each $n\in \Z$, there exists an automorphism
 $\rho_{n}$ of $V_{\hat{B}_{q}}(\ell,0)$, which is uniquely determined by
\begin{eqnarray}
\rho_{n}(\beta^{(r)})=\beta^{(n+r)}\ \ \ \mbox{ for }r\in \Z.
\end{eqnarray}
Furthermore, $\rho_{0}=1$ and $\rho_{m+n}=\rho_{m}\rho_{n}$ for
$m,n\in \Z$.
 \el

\begin{proof} First, for each given $n\in \Z$, from the defining relations
of $\hat{B}_{q}$ we see that there exists an automorphism
$\sigma_{n}$ of $\hat{B}_{q}$, which is uniquely determined by
\begin{eqnarray*}
\sigma_{n}(c)=c,\ \ \ \
\sigma_{n}(\beta^{(r)}_{m})=\beta^{(n+r)}_{m} \ \ \ \mbox{ for
}r,m\in \Z.
\end{eqnarray*}
Then $\sigma_{n}$ induces an automorphism of $U(\hat{B}_{q})$.
Clearly, $\sigma_{n}(\hat{B}_{q}^{+})=\hat{B}_{q}^{+}$. It follows
that $\sigma_{n}$ reduces to a linear automorphism $\rho_{n}$ of
$V_{\hat{B}_{q}}(\ell,0)$ with $\rho_{n}({\bf 1})={\bf 1}$ and
$$\rho_{n}(av)=\sigma_{n}(a)\rho_{n}(v)\ \ \mbox{ for }a\in U(\hat{B}_{q}),\ v\in
V_{\hat{B}_{q}}(\ell,0).$$ For $n,r\in \Z$, we have
$$\rho_{n}(\beta^{(r)})=\rho_{n}(\beta^{(r)}_{-1}{\bf
1})=\sigma_{n}(\beta^{(r)}_{-1})\rho_{n}({\bf 1})
=\beta^{(n+r)}_{-1}{\bf 1}=\beta^{(n+r)}.$$
$$\rho_{n}(Y(\beta^{(r)},x)v)=\rho_{n}(\beta^{(r)}(x)v)=\beta^{(n+r)}(x)\rho_{n}(v)=
Y(\rho_{n}(\beta^{(r)}),x)\rho_{n}(v)$$ for $v\in
V_{\hat{B}_{q}}(\ell,0)$. Consequently, $\rho_{n}$ is an
automorphism of weak quantum vertex algebra
$V_{\hat{B}_{q}}(\ell,0)$.
\end{proof}

Next, we relate restricted $\widetilde{\mathcal{H}}_{q}$-modules of
level $\ell$ to $\phi$-coordinated quasi
$V_{\hat{B}_{q}}(\ell,0)$-modules. To achieve this goal we shall
need the following result:

\bl{lspecial} Let $W$ be a vector space and let $K$ be any
$Y_{\E}^{e}$-closed quasi compatible subspace of $\E(W)$ with
$a(x),b(x)\in K$. Suppose that $a(x),b(x)$ satisfy relation
\begin{eqnarray}
[a(x_{1}),b(x_{2})]=f(x_{1}/x_{2})+\alpha
\left(x_{2}\frac{\partial}{\partial x_{2}}\right)
\delta\left(\frac{x_{2}}{x_{1}}\right)+\beta\left(x_{2}\frac{\partial}{\partial
x_{2}}\right)\delta\left(\frac{qx_{2}}{x_{1}}\right),
\end{eqnarray}
where $f(x)\in \C(x),\; \alpha,\beta, q\in \C$ with $q\ne 0,1$. Then
\begin{eqnarray}
[Y_{\E}^{e}(a(x),x_{1}),Y_{\E}^{e}(b(x),x_{2})]
=\iota_{x_{2},x_{1}}(f(e^{x_{1}-x_{2}}))+\alpha
\frac{\partial}{\partial x_{2}}x_{1}^{-1}
\delta\left(\frac{x_{2}}{x_{1}}\right).
\end{eqnarray}
\el

\begin{proof} {}From the given relation we have
\begin{eqnarray}\label{equasi-locality}
&&(x_{1}-x_{2})^{2}(x_{1}-qx_{2})^{2}a(x_{1})b(x_{2})\nonumber\\
&=&(x_{1}-x_{2})^{2}(x_{1}-qx_{2})^{2}\left(b(x_{2})a(x_{1})+f(x_{1}/x_{2})\right),
\end{eqnarray}
which implies
$$(x_{1}-x_{2})^{2}(x_{1}-qx_{2})^{2}a(x_{1})b(x_{2})\in
\Hom(W,W((x_{1},x_{2}))).$$ Then we have
\begin{eqnarray*}
&&x^{4}(e^{z}-1)^{2}(e^{z}-q)^{2}Y_{\E}^{e}(a(x),z)b(x)\nonumber\\
&=&\left((x_{1}-x)^{2}(x_{1}-qx)^{2}a(x_{1})b(x)\right)|_{x_{1}=xe^{z}}\nonumber\\
&=&\Res_{x_{1}}x_{1}^{-1}\delta\left(\frac{xe^{z}}{x_{1}}\right)
\left((x_{1}-x)^{2}(x_{1}-qx)^{2}a(x_{1})b(x)\right)\nonumber\\
&=&\Res_{x_{1}}\frac{1}{x_{1}-xe^{z}}\left((x_{1}-x)^{2}(x_{1}-qx)^{2}a(x_{1})b(x)\right)\nonumber\\
&&-\Res_{x_{1}}\frac{1}{-xe^{z}+x_{1}}\left((x_{1}-x)^{2}(x_{1}-qx)^{2}a(x_{1})b(x)\right)\nonumber\\
&=&\Res_{x_{1}}\frac{1}{x_{1}-xe^{z}}(x_{1}-x)^{2}(x_{1}-qx)^{2}a(x_{1})b(x)\nonumber\\
&&-\Res_{x_{1}}\frac{1}{-xe^{z}+x_{1}}(x_{1}-x)^{2}(x_{1}-qx)^{2}(b(x)a(x_{1})+f(x_{1}/x)).
\end{eqnarray*}
For convenience let us denote the last term by $A$. Noticing that
$A$ involves only nonnegative powers of $z$, we have
\begin{eqnarray}\label{eyabz}
Y_{\E}^{e}(a(x),z)b(x)
&=&x^{-4}\iota_{z,0}\left((e^{z}-1)^{-2}(e^{z}-q)^{-2}\right)A\nonumber\\
&=&x^{-4}\left(\frac{1}{(1-q)^{2}}z^{-2}+\frac{q-3}{(1-q)^{3}}z^{-1}\right)A+O(z^{0}).
\end{eqnarray}
As
\begin{eqnarray*}
&&\frac{1}{x_{1}-xe^{z}}(x_{1}-x)^{2}=(x_{1}-x)+xz+O(z^{2}),\\
&&\frac{1}{-xe^{z}+x_{1}}(x_{1}-x)^{2}=(x_{1}-x)+xz+O(z^{2}),\\
&&(x_{1}-qx_{2})^{2}\left(x_{2}\frac{\partial}{\partial
x_{2}}\right)\delta\left(\frac{qx_{2}}{x_{1}}\right)=0,
\end{eqnarray*}
we get
\begin{eqnarray*}
A&=&\Res_{x_{1}}(x_{1}-x+xz)(x_{1}-qx)^{2}
\left(a(x_{1})b(x)-b(x)a(x_{1})-f(x_{1}/x)\right)+O(z^{2})\nonumber\\
&=&\alpha \Res_{x_{1}}(x_{1}-x+xz)(x_{1}-qx)^{2}
\left(x\frac{\partial}{\partial x}\right)
\delta\left(\frac{x}{x_{1}}\right)+O(z^{2})\\
&=&\alpha x^{4} \left( (1-q)^{2}+z(1-q)(3-q)\right)+O(z^{2}).
\end{eqnarray*}
Combining this with (\ref{eyabz}) we obtain
\begin{eqnarray*}
Y_{\E}^{e}(a(x),z)b(x)=\alpha z^{-2}+O(z^{0}),
\end{eqnarray*}
which implies
\begin{eqnarray}\label{esingular-relation}
a(x)_{n}^{e}b(x)=\delta_{n,1}\alpha 1_{W} \ \ \ \mbox{ for }n\ge 0.
\end{eqnarray}
With (\ref{equasi-locality}), by Proposition 5.3 of
\cite{li-phi-module} (specially by (5.12) in the proof) we have
\begin{eqnarray*}
&&Y_{\E}^{e}(a(x),x_{1})Y_{\E}^{e}(b(x),x_{2})-Y_{\E}^{e}(b(x),x_{2})Y_{\E}^{e}(a(x),x_{1})
-\iota_{x_{2},x_{1}}f(e^{x_{1}-x_{2}})\nonumber\\
&=&\sum_{n\ge
0}Y_{\E}^{e}(a(x)_{n}^{e}b(x),x_{2})\frac{1}{n!}\left(\frac{\partial}{\partial
x_{2}}\right)^{n}x_{1}^{-1}\delta\left(\frac{x_{2}}{x_{1}}\right).
\end{eqnarray*}
Then the desired relation follows from (\ref{esingular-relation}).
\end{proof}

 Recall $\Gamma_{q}=\{ q^{n}\;|\;
n\in \Z\}\subset \C^{\times}$.

\bt{theisenberg} Let $W$ be a restricted
$\widetilde{\mathcal{H}}_{q}$-module of level $\ell\in \C$. There
exists a $\Gamma_{q}$-covariant $\phi$-coordinated quasi module
structure on $W$ for $V_{\hat{B}_{q}}(\ell,0)$, which is uniquely
determined by the condition that
$Y_{W}(\beta^{(r)},x)=\beta(q^{r}x)$ for $r\in \Z$. \et

\begin{proof} Set
$$U_{W}=\{1_{W}\}\cup \{ \beta(q^{r}x)\;|\; r\in \Z\}\subset \E(W).$$
Let $a,b\in \Gamma_{q}$. {}From (\ref{eqlocality-commutator}) we
have
\begin{eqnarray}\label{ebetaabrelation}
&&\beta(ax_{1})\beta(bx_{2})-\beta(bx_{2})\beta(ax_{1})\nonumber\\
&=&\ell\left(\frac{qab^{-1}x_{1}/x_{2}}{(1-qab^{-1}x_{1}/x_{2})^{2}}
-\frac{q^{-1}ab^{-1}x_{1}/x_{2}}{(1-q^{-1}ab^{-1}x_{1}/x_{2})^{2}}
\right)\ \ \nonumber\\
&&+\ell \left(x_{2}\frac{\partial}{\partial x_{2}}\right)\left(
\delta\left(\frac{ba^{-1}x_{2}}{x_{1}}\right)
-\delta\left(\frac{qba^{-1}x_{2}}{x_{1}}\right)\right).
\end{eqnarray}
{}From this we see that $U_{W}$ is a quasi $\S_{trig}$-local subset
of $\E(W)$. By Theorem \ref{trecall1}, $U_{W}$ generates a weak
quantum vertex algebra $\<U_{W}\>_{e}$ with $W$ as a faithful
$\phi$-coordinated quasi module.

With (\ref{ebetaabrelation}), by Lemma \ref{lspecial} we have
\begin{eqnarray}
&&Y_{\E}^{e}(\beta(ax),x_{1})Y_{\E}^{e}(\beta(bx),x_{2})
-Y_{\E}^{e}(\beta(bx),x_{2})Y_{\E}^{e}(\beta(ax),x_{1})\nonumber\\
&=&
\iota_{x_{2},x_{1}}\left(\frac{qab^{-1}e^{x_{1}-x_{2}}}{(1-qab^{-1}e^{x_{1}-x_{2}})^{2}}
-\frac{q^{-1}ab^{-1}e^{x_{1}-x_{2}}}{(1-q^{-1}ab^{-1}e^{x_{1}-x_{2}})^{2}}
\right)\ell\nonumber\\
&&+\ell (\delta_{a,b}-\delta_{a,qb})\left(\frac{\partial}{\partial
x_{2}}\right) x_{1}^{-1}\delta\left(\frac{x_{2}}{x_{1}}\right).
\end{eqnarray}
It follows that $\<U_{W}\>_{e}$ is a restricted $\hat{B}_{q}$-module
of level $\ell$ with $\beta^{(r)}(z)$ acting as
$Y_{\E}^{e}(\beta(q^{r}x),z)$ for $r\in \Z$. We also have
$Y_{\E}^{e}(\beta(q^{r}x),z)1_{W}\in \<U_{W}\>_{e}[[z]]$ for $r\in
\Z$. Then $\<U_{W}\>_{e}$ together with $1_{W}$ is a vacuum
$\hat{B}_{q}$-module of level $\ell$. By the universality of
$V_{\hat{B}_{q}}(\ell,0)$, there exists a $\hat{B}_{q}$-module
homomorphism $\psi$ from $V_{\hat{B}_{q}}(\ell,0)$ to
$\<U_{W}\>_{e}$ with $\psi({\bf 1})=1_{W}$. For $r\in \Z$, we have
$$\psi(\beta^{(r)})=\psi(\beta^{(r)}_{-1}{\bf 1})=\beta(q^{r}x)^{e}_{-1}1_{W}
=\beta(q^{r}x).$$
 Furthermore, we have
$$\psi(Y(\beta^{(r)},z)v)=\psi(\beta^{(r)}(z)v)=Y_{\E}^{e}(\beta(q^{r}x),z)\psi(v)
=Y_{\E}^{e}(\psi(\beta^{(r)}),z)\psi(v)$$ for $r\in \Z,\; v\in
V_{\hat{B}_{q}}(\ell,0)$.
 Since
$V_{\hat{B}_{q}}(\ell,0)$ as a weak quantum vertex algebra is
generated by elements $\beta^{(r)}$ $(r\in \Z)$, it follows that
$\psi$ is a homomorphism of  weak quantum vertex algebras.
Consequently, $W$ is a $\phi$-coordinated quasi
$V_{\hat{B}_{q}}(\ell,0)$-module where for $r\in \Z$,
$$Y_{W}(\beta^{(r)},z)=Y_{W}(\psi(\beta^{(r)}),z)=Y_{W}(\beta(q^{r}x),z)=\beta(q^{r}z).$$

For $n,r\in \Z$, we have
$$Y_{W}(\rho_{n}(\beta^{(r)}),x)=Y_{W}(\beta^{(n+r)},x)
=\beta(q^{n+r}x)=Y_{W}(\beta^{(r)},q^{n}x).$$
Since
$V_{\hat{B}_{q}}(\ell,0)$ is generated by elements $\beta^{(r)}$ for
$r\in \Z$, the $\Gamma_{q}$-covariance follows from Lemma
\ref{lsimplefacts}. Therefore, $W$ is a $\Gamma_{q}$-covariant
$\phi$-coordinated quasi $V_{\hat{B}_{q}}(\ell,0)$-module.
\end{proof}

\end{document}